\documentclass{amsart}
\usepackage{amsfonts}
\usepackage{amsmath,amssymb}
\usepackage{amsthm}
\usepackage{amscd}
\usepackage{graphics}
\usepackage{graphicx}

\theoremstyle{remark}{
\newtheorem{Def}{{\rm Definition}}
\newtheorem{Ex}{{\rm Example}}

\newtheorem{Prob}{{\rm Problem}}

}
\theoremstyle{plain}
{

\newtheorem{MainThm}{Main Theorem}

}

\begin{document}
\title[Refined algebraic domains with differential geometric finite sets]{Refined algebraic domains with finite sets in the boundaries respecting differential geometry}
\author{Naoki kitazawa}

\address{
}
\email{naokikitazawa.formath@gmail.com}
\urladdr{https://naokikitazawa.github.io/NaokiKitazawa.html}
\maketitle
\begin{abstract}
We are interested in shapes of real algebraic curves in the plane and regions surrounded by them: they are named {\it refined algebraic domains} by the author.

As characteristic finite sets, we consider points contained in two curves and the sets of singular points of the restrictions of the projections to the lines to the curves. As a new case, we respect differential geometry and consider inflection points and points of some double tangent lines of a single connected curve. We prove fundamental properties and investigate some examples. 

We have also previously considered the cases where the curves are straight lines, circles, or boundaries of ellipsoids for example. Such simple cases are trivial in our new consideration.

\end{abstract}
%【REVISE】 combinatoric ～ is → combinatorial object. It is .
%【REVISE】  such that a point is a vertex if and only if the corresponding connected component of the level set contains some singular points → whose vertex set is the set of all points containing some singular points in the corresponding connected component of the level set .
%【REVISE】 We delete "extending the result before".
\section{Introduction.}
\label{sec:1}
Our interest lies in arrangements of real algebraic curves in the plane and regions surrounded by them and their shapes ({\it refined algebraic domains}: \cite{kitazawa6}). We are also interested in characteristic finite sets in the boundaries of the regions. 

This is natural in mathematics, especially in real algebraic geometry and combinatorics, for example. Such studies are also new, surprisingly.
(Some of) related studies have been started by the author. As a natural question, we have been interested in explicit construction of explicit objects in real algebraic geometry: real algebraic manifolds and maps, especially ones generalizing natural projections of spheres. 
We have considered and succeeded in reconstruction of real algebraic maps onto given regions surrounded by real algebraic curves \cite{kitazawa1, kitazawa2}. This is also different from approximations by or existence of real algebraic objects: positive classical results are well-known, thanks to Weierstrass, followed by fundamental classical studies in real algebraic geometry, by Nash and Tognoli, for example.
\subsection{Fundamental notation and notions.}
Let ${\mathbb{R}}^k$ denote the $k$-dimensional Euclidean space. This is a smooth manifold and a Riemannian manifold equipped with the standard Euclidean metric. Let $||x||\geq 0$ denote the distance of $x$ and the origin $0$ under the metric. We regard the space as the k-dimensional real vector space canonically. Let ${\mathbb{R}}^1;=\mathbb{R}$. 
Let ${\pi}_{m,n}:{\mathbb{R}}^m \rightarrow {\mathbb{R}}^n$ ($m>n \geq 1$) denote the so-called canonical projection ${\pi}_{m,n}(x)=x_1$ ($x=(x_1,x_2) \in {\mathbb{R}}^n \times {\mathbb{R}}^{m-n}={\mathbb{R}}^m$). 
Let $D^k:=\{x \in {\mathbb{R}}^k \mid ||x|| \leq 1\}$ and $S^k:=\{x \in {\mathbb{R}}^{k+1} \mid ||x||=1\}$: they are the $k$-dimensional unit disk and the $k$-dimensional unit sphere. 

Let $X$ be a topological space and $Y \subset X$ a subspace. We use $\overline{Y}$ for the closure and $Y^{\circ}$ for the interior in $X$. In our paper, $X$ is the Euclidean space (mainly ${\mathbb{R}}^2$) considered in the discussions unless otherwise stated. For a topological space $X$ regarded as a so-called cell complex, the dimension $\dim X$ can be defined uniquely as the dimension of the cell of the maximal dimension. This only depends on the topology of $X$. Topological manifolds, polyhedra, and graphs, or more generally, $1$-dimensional CW complexes, are of such spaces. %For a topological manifold $X$ whose boundary is non-empty, we use $\partial X$ for the boundary and ${\rm Int}\ X:=X-\partial X$. 
For a smooth manifold $X$ and a point $x \in X$, let $T_xX$ denote the tangent vector space of $X$ at $x$. For smooth manifolds $X$ and $Y$ and a smooth map $c:X \rightarrow Y$, $x \in X$ is a {\it singular point} of $c$ if the rank of the differential ${dc}_x:T_xX \rightarrow T_{c(x)}Y$ is smaller than both the dimensions $\dim X$ and $\dim Y$: note that the differential ${dc}_x$ at $x \in X$ is linear.
A union $S$ of connected components of the zero set of a real polynomial map is {\it non-singular} if the polynomial map has no singular point in $S$. This respects the implicit function theorem. 
A {\it straight line} is the zero set of a real polynomial of degree $1$ and non-singular and a straight segment is a subset there diffeomorphic to $D^1$. A {\it circle} $S_{x_0,r}:=\{x \in {\mathbb{R}}^2 \mid ||x-x_0||=r\}$ ($r>0$) and the boundary of an {\it ellipsoid} $S_{x_0,r_1,r_2}:=\{x=(x_1,x_2) \in {\mathbb{R}}^{2} \mid \frac{{(x_1-x_{0,1})}^2}{r_1}+\frac{{(x_2-x_{0,2})}^2}{r_2}=1\}$ ($x_0=(x_{0,1},x_{0,2}) \in {\mathbb{R}}^2$, $r_j>0$) are the zero sets of real polynomials of degree $2$ and non-singular. Remember $S^1=S_{0,1}$ for example.

A straight line $L \subset {\mathbb{R}}^2$ is {\it tangent to }{\rm (}{\it the image of}{\rm )}{\it a smooth embedding $S \subset {\mathbb{R}}^2$ of a $1$-dimensional manifold (or a smooth curve $S$) at $p \in S \bigcap L$} or the {\it tangent line of $S$ at $p$} if the tangent vector spaces of $L$ and $S$ at $p$ coincide. The point $p \in S$ is also a {\it reflection point of $S$} if $L$ is tangent to $S$ at $p$ and for any open neighborhood $U_p \subset {\mathbb{R}}^2$ of $p$, $(L-\{p\}) \bigcap U_p$ is not contained in a single connected component of $U_p-(S \bigcap U_p)$ or in the curve $S$. In other words, the curve $S$ is locally the graph of a smooth function such that $p$ is a singular point of it and that the 2nd derivative of the function is $0$ at $p$. 
Note also that the tangent line of $S$ at $p$ is in the curve if and only if $S$ is a straight line. A {\it double tangent line of a connected smooth curve $S$} is a straight line tangent to $S$ at least two points of $S$.
Straight lines, circles, and the boundary of an ellipsoid have no reflection points or double tangent lines. 

A {\it graph} is a CW complex of dimension $1$ with $1$-cells ({\it edges}) and $0$-cells ({\it vertices}). The set of all edges (vertices) of the graph is the {\it edge set} ({\it vertex set}) of it. Two graphs $G_1$ and $G_2$ are {\it isomorphic} if there exists a (piecewise smooth) homeomorphism $\phi:G_1 \rightarrow G_2$ mapping the vertex set of $G_1$ onto that of $G_2$ and this is called an {\it isomorphism} of the graphs. A {\it digraph} is a graph all of whose edges are oriented. Two digraphs are {\it isomorphic} if there exists an isomorphism of graphs preserving the orientations and this is called an {\it isomorphism} of the digraphs. For a digraph $G$ and a map $V_{G}$ on its vertex set onto a partially ordered set $P$ with natural conditions, we can orient the graph according to the values: each edge $e$ of the graph connects two distinct vertices $v_{e,1}$ and $v_{e,2}$ and is oriented according to the rule that the edge $e$ departs from $v_{e,1}$ and enters $v_{e,2}$ if and only if $V_{G}(v_{e,1})<V_{G}(v_{e,2})$ where "$<$" denotes the order on $P$. 
A pair of this graph $G$ and a map $V_{G}$ is  {\it a V-digraph}. For V-graphs, {\it isomorphisms} between two V-digraphs with the relation that two V-digraphs are {\it isomorphic} are canonically defined.
\subsection{Refined algebraic domains, their characteristic finite sets from differential geometry, and our main work.}
\label{subsec:1.2}
%We introduce {\it arrangements of circles for Morse-Bott functions} or {\it MBC arrangements}
The following has been defined first in \cite{kitazawa6} and also used in \cite{kitazawa7} where some small difference exists between the definitions.
\begin{Def}
\label{def:1}
A pair of a family $\mathcal{S}=\{S_j \subset {\mathbb{R}}^2\}$ of the zero set of a real polynomial being non-singular and connected and a region $D_{\mathcal{S}} \subset {\mathbb{R}}^2$ satisfying the following is a {\it refined algebraic domain}. 
\begin{enumerate}
\item \label{def:1.1} The region $D_{\mathcal{S}}$ is a connected component of ${\mathbb{R}}^2-{\bigcup}_{S_j \in \mathcal{S}} S_j$ such
 that the intersection $\overline{D_{\mathcal{S}}} \bigcap S_j$ is not empty for any curve $S_j \in \mathcal{S}$.
\item \label{def:1.2} At each point $p_{j_1,j_2} \in \overline{D_{\mathcal{S}}}$, at most two distinct curves $S_{j_1}, S_{j_2} \in \mathcal{S}$ intersect as follows (in other words transversality is satisfied).
For each point $p_{j_1,j_2}$ before, the sum of the tangent vector spaces of them at $p_{j_1,j_2}$ and the tangent vector space of ${\mathbb{R}}^2$ at $p_{j_1,j_2}$ agree. %and three distinct curves do not intersect in $\overline{D_{\mathcal{S}}}$.
\end{enumerate}
\end{Def}
For example, we first define the characteristic finite set $F_{D_{\mathcal{S}},1}:=F_{D_{\mathcal{S}},0} \sqcup F_{D_{\mathcal{S}},1,{\rm pole}}$ from singularity theory as follows.

\begin{itemize}
\item $F_{D_{\mathcal{S}},0}$ is the set of all points in $\overline{D_{\mathcal{S}}}$ in two distinct curves $S_{j_1}, S_{j_2} \in \mathcal{S}$.
\item By removing the set of all points above from the set $\overline{D_{\mathcal{S}}}-D_{\mathcal{S}}$, we have a smooth manifold of dimension $1$ or a smooth curve with no boundary. We define $F_{D_{\mathcal{S}},1,{\rm pole}}$ to be the set of all points which are also singular points of the restriction of ${\pi}_{2,1}$ to the obtained smooth curve in $\overline{D_{\mathcal{S}}}-D_{\mathcal{S}}$ and which are also isolated in the obtained smooth curve. 
\end{itemize}
Note that in this situation if $p \in \overline{D_{\mathcal{S}}}-D_{\mathcal{S}}$ is in a curve $S_j \in \mathcal{S}$ of the form $S_{t}:=\{x=(x_1,x_2) \mid x_1=t \}$, then this is a singular point of the restriction of ${\pi}_{2,1}$ to the obtained smooth curve in $\overline{D_{\mathcal{S}}}-D_{\mathcal{S}}$ and this is not in $F_{D_{\mathcal{S}},1,{\rm pole}}$.
 
We have the equivalence relation ${\sim}_{D_{\mathcal{S}},1}$ on $\overline{D_{\mathcal{S}}}$: two points are equivalent if and only if they are in a same component of the preimage of a same point for the restriction of ${\pi}_{2,1}$ to $\overline{D_{\mathcal{S}}}$. Let $q_{D_{\mathcal{S}},1}$ denote the quotient map. We have the function $V_{D_{\mathcal{S}},1}$ satisfying ${\pi}_{2,1}=V_{D_{\mathcal{S}},1} \circ q_{D_{\mathcal{S}},1}$ uniquely. The quotient space $W_{D_{\mathcal{S}},1}:=\overline{D_{\mathcal{S}}}/{\sim}_{D_{\mathcal{S}},1}$ is a V-digraph. We can check this from \cite{saeki1, saeki2} or see \cite{kitazawa4, kitazawa5} for example. We can also check more explicitly. We do not need to understand the theory.
\begin{enumerate}
\item The vertex set is the set of all points $v$ such that the preimage ${q_{D_{\mathcal{S}},1}}^{-1}(v)$ contains some point of $F_{D_{\mathcal{S}},1}$.
\item The edge connecting $v_1$ and $v_2$ are oriented as one departing from $v_1$ and entering $v_2$ according to $V_{D_{\mathcal{S}},1}(v_1)<V_{D_{\mathcal{S}},1}(v_2)$.
\end{enumerate}
\begin{Def}
\label{def:2}
We call the (V-di)graph $(W_{D_{\mathcal{S}},1},V_{D_{\mathcal{S}},1})$ a {\it Poincar\'e-Reeb} ({\it V-di}){\it graph of $D_{\mathcal{S}}$}. We omit the function $V_{D_{\mathcal{S}},1}$ if we can guess easily.
\end{Def}
Related original studies are shown in \cite{bodinpopescupampusorea, sorea1, sorea2}.
In our paper, we consider the characteristic finite sets from differential geometry of curves first.
\begin{Def}
If for a refined algebraic domain $(\mathcal{S},D_{\mathcal{S}})$, no reflection point of a curve $S_j \in \mathcal{S}$ being also a point of $\overline{D_{\mathcal{S}}}-D_{\mathcal{S}}$ exists for any $S_j \in \mathcal{S}$, then the refined algebraic domain is {\it with no inflection point} or {\it NIP}. If for a refined algebraic domain $(\mathcal{S},D_{\mathcal{S}})$, no point of a curve $S_j \in \mathcal{S}$ being also contained in some double tangent line of a curve $S_j \in \mathcal{S}$ and being also a point of $\overline{D_{\mathcal{S}}}-D_{\mathcal{S}}$ exists for any $S_j \in \mathcal{S}$, then the refined algebraic domain is {\it with no double tangent lines} or {\it NDTL}.
\end{Def}
The next section explains our new result. Main Theorem \ref{mthm:1} is for the existence of an NIP refined algebraic domain whose Poincar\'e-Reeb graph is isomorphic to a given graph of a certain class. This extends or is regarded as a variant of our previous result \cite{kitazawa6}. For NDTL, we have Main Theorem \ref{mthm:2}.

%Hereafter, the Euclidean space ${\mathbb{R}}^k$ is regarded as the vector space canonically.

\section{Our new result.}
\label{sec:2}
Hereafter, we use ${\pi}_{m,n,i}:{\mathbb{R}}^m \rightarrow {\mathbb{R}}^n$ ($m>n$) for the projection defined by ${\pi}_{m,n,i}(x):=x_i$ ($x=(x_1,x_2) \in {\mathbb{R}}^n \times {\mathbb{R}}^{m-n}={\mathbb{R}}^m$) for $i=1,2$. Remember ${\pi}_{m,n}={\pi}_{m,n,1}$.
We can define the characteristic finite set $F_{D_{\mathcal{S}},2}:=F_{D_{\mathcal{S}},0} \sqcup F_{D_{\mathcal{S}},2,{\rm pole}}$ in the same way for ${\pi}_{2,1,2}$. We can also define similar objects and we abuse the natural notation.
We can have the (V-di)graph $(W_{D_{\mathcal{S}},2},V_{D_{\mathcal{S}},2})$ a {\it Poincar\'e-Reeb }({\it V-di}){\it graph of $D_{\mathcal{S}}$ for ${\pi}_{2,1,2}$}. 
We omit the function $V_{D_{\mathcal{S}},2}$ if we can guess easily.
We can also call the V-(di)graph $(W_{D_{\mathcal{S}},1},V_{D_{\mathcal{S}},1})$ the {\it Poincar\'e-Reeb }({\it V-di}){\it graph of $D_{\mathcal{S}}$ for ${\pi}_{2,1,1}$}.
\begin{Def}
A refined algebraic domain $(\mathcal{S},D_{\mathcal{S}})$ is said to be {\it Morse} if at each point of $F_{D_{\mathcal{S}},i,{\rm pole}}$ the value of the 2nd derivative of the restriction of ${\pi}_{2,1,i}$ to the smooth curve $\overline{D_{\mathcal{S}}}-D_{\mathcal{S}}$ is not $0$ for $i=1,2$ and at each point $p$ of $F_{D_{\mathcal{S}},0}$ the values of the 1st derivative of the restriction of ${\pi}_{2,1,i}$ to the two smooth curves $S_{j_1}$ and $S_{j_2}$ containing the point $p$ are not $0$ for $i=1,2$.
A Morse refined algebraic domain is also {\it generic} or {\it G-Morse} if the restriction of ${\pi}_{2,1,i}$ to  $F_{D_{\mathcal{S}},i}$ is injective for $i=1,2$.
\end{Def}
\begin{MainThm}
\label{mthm:1}
Let a refined algebraic domain $(\mathcal{S},D_{\mathcal{S}})$ be Morse. Then we can have another Morse refined algebraic domain $({\mathcal{S}}^{\prime},D_{{\mathcal{S}}^{\prime}})$ enjoying the following.
\begin{enumerate}
\item \label{mthm:1.1} The new refined algebraic domain $({\mathcal{S}}^{\prime},D_{{\mathcal{S}}^{\prime}})$ is NIP.
\item \label{mthm:1.2} The relations $\mathcal{S} \subset {\mathcal{S}}^{\prime}$ and $D_{{\mathcal{S}}^{\prime}} \subset D_{\mathcal{S}}$ hold.
\item \label{mthm:1.3} The V-digraphs $(W_{D_{\mathcal{S}},i},V_{D_{\mathcal{S}},i})$ and $(W_{D_{{\mathcal{S}}^{\prime}},i},V_{D_{{\mathcal{S}}^{prime}},i})$ are homeomorphic for $i=1,2$.
\end{enumerate}

\end{MainThm}
Hereafter, we apply approximation of smooth curves by real algebraic curves. For real algebraic geometry, see \cite{bochnakcosteroy, kollar}. For the usage of approximation, see also \cite{bodinpopescupampusorea, elredge, lellis}.
\begin{proof}[A proof of Main Theorem \ref{mthm:1}]
Each inflection point $p$ of each curve $S_j$ must not be points of $F_{D_{\mathcal{S}},1,{\rm pole}} \bigcup F_{D_{\mathcal{S}},2,{\rm pole}}$.

We choose an inflection point $p$ of a curve $S_j \in \mathcal{S}$ such that $p$ is also in the closure $\overline{D_{\mathcal{S}}}$.
We choose such an inflection point of $S_j$ contained in exactly one curve $S_j \in \mathcal{S}$. 
We apply new technique respecting \cite{kitazawa5} where we do not need to understand the study. For each inflection point $p$ of the curve here, we consider a small ellipsoid containing the point in the interior and a disk $D_p$ obtained by approximation of the ellipsoid respecting the derivations of the class $k \geq 2$. We can also choose the disk $D_p$ in such a way that the intersection of the boundary $\partial D_p$ of the new disk and the region $\overline{D_{\mathcal{S}}}$ is represented as a smooth connected curve sufficiently close to the point $p$ and containing two points $p_{\rm l}, p_{\rm r} \in \overline{D_{\mathcal{S}}}-D_{\mathcal{S}}$ different from and sufficiently close to $p$. We can also do this in such a way that the boundary of the new disk $D_p$ is also the zero set of a real polynomial function and non-singular and that the condition (\ref{def:1.2}) of Definition \ref{def:1} or the condition on the intersections of the curves surrounding the region is satisfied. See also FIGURE \ref{fig:1}. For example, the green colored curve connecting the two green colored points $p_{\rm l} \neq p$ and $p_{\rm r} \neq p$ in the first figure is the curve in the closed curve $\partial D_p$.
\begin{figure}
\includegraphics[height=75mm, width=100mm]{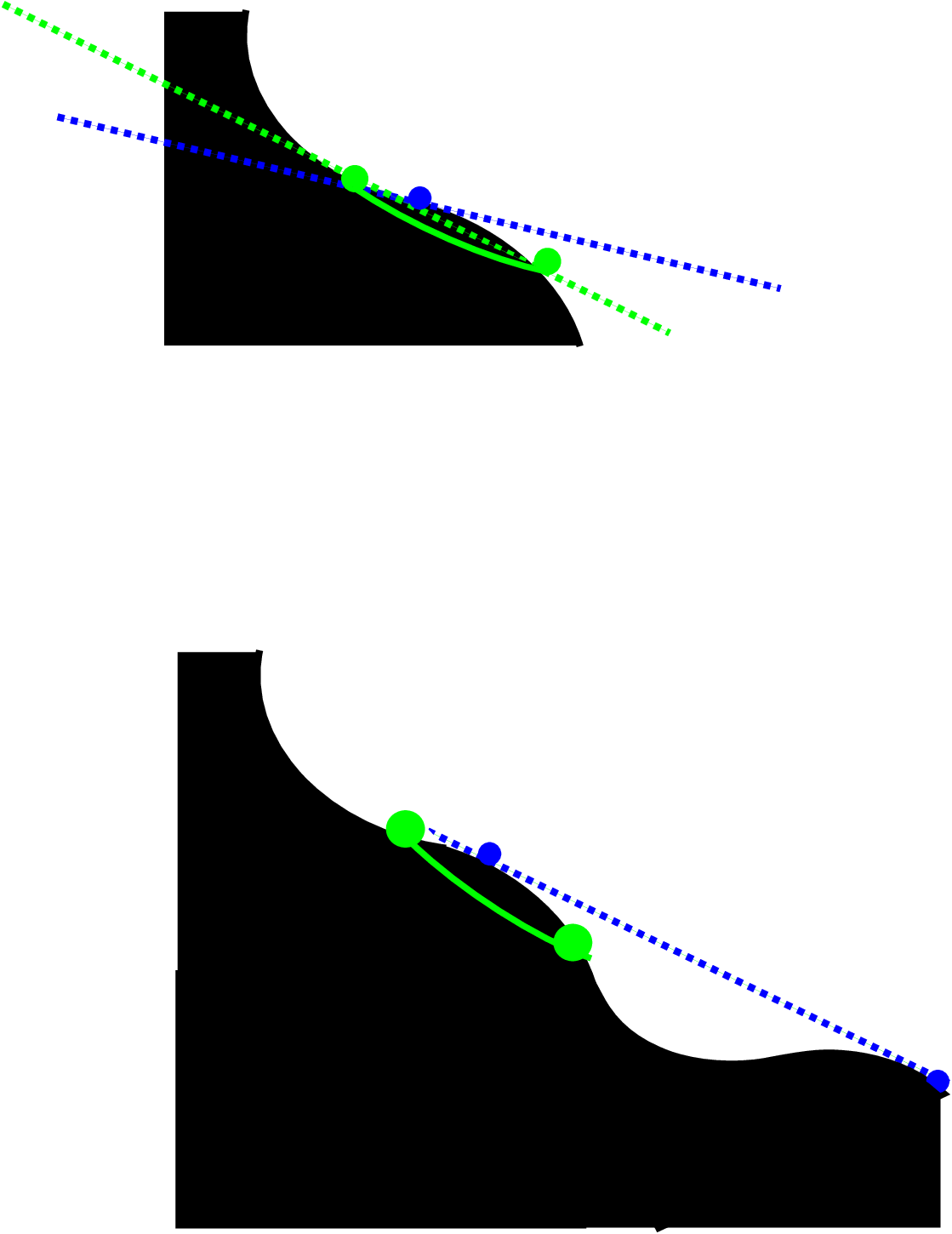}
\caption{The region $\overline{D_{\mathcal{S}}}$ is colored in black. The blue colored point shows an inflection point $p$ of a smooth connected curve $S_j$ in the first case. The blue colored point shows a point in a double tangent line of the smooth connected curve $S_j$ in the second case. These points $p$ are also points of $\overline{D_{\mathcal{S}}}-D_{\mathcal{S}}$. 
Blue colored dotted segments show tangent lines to the points $p$.
The two points $p_{\rm l} \neq p$ and $p_{\rm r} \neq p$ are colored in green. In the first case, they are connected by a smooth curve, colored in green, and sufficiently close to $p$.}
\label{fig:1}
\end{figure}

Note again that the boundary of an ellipsoid does not contain inflection points. Note also that our approximation does not change this property.

We choose such an inflection point of $S_j$ contained in $F_{D_{\mathcal{S}},0}$. We can apply similar technique and discuss similarly. See also FIGURE \ref{fig:2}. We can also abuse the same notation again. For example, we can choose an ellipsoid, consider its suitable approximation to have another disk $D_p$, and find a curve which connects the two green colored points $p_{\rm l} \neq p$ and $p_{\rm r} \neq p$ and which is in the closed curve $\partial D_p$ as before. Note also that we can and we also need to choose the two points in such a way that the straight segment connecting them is not in the straight line of the form $\{(t,y)\mid t \in \mathbb{R}\}$ or $\{(y,t)\mid t \in \mathbb{R}\}$. If we need, then refer to the original preprint \cite{kitazawa5}.
By adding the disks, we have our desired case, enjoying (\ref{mthm:1.1})--(\ref{mthm:1.3}). This also includes a kind of routine work. Especially, for the property (\ref{mthm:1.3}), it is important that adding the disks and removing some subsets from the given regions does not change the local topologies of our graphs. If we need, then refer to the original study (\cite{kitazawa5}) again.
 
\begin{figure}
\includegraphics[height=75mm, width=100mm]{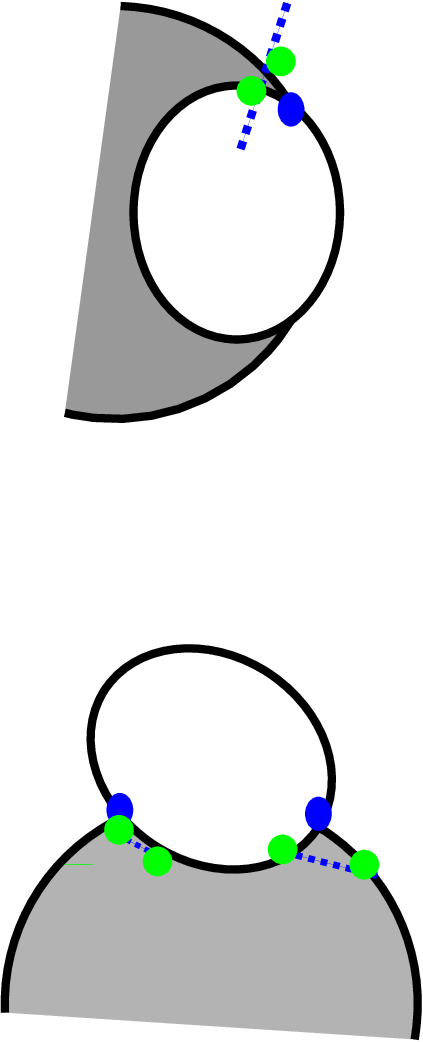}
\caption{The regions $\overline{D_{\mathcal{S}}}$ are colored in gray. The blue colored points show the points of $p \in F_{D_{\mathcal{S}},0}$ which may be inflection points of smooth connected curves $S_j$ or points in some double tangent lines of the smooth connected curves $S_j$. The two points $p_{\rm l} \neq p$ and $p_{\rm r} \neq p$ are colored in green and connected by straight segments, dotted in blue and sufficiently close to $p$.}
	\label{fig:2}
\end{figure}

This completes the proof.
\end{proof}
\begin{Ex}[\cite{bodinpopescupampusorea, kitazawa6}]
\label{ex:1}
Let $c_G:G \rightarrow \mathbb{R}$ be a piecewise smooth function on a finite and connected graph $G$ with the following.
\begin{enumerate}
\item The restriction of $c_G$ to each edge is injective.
\item A vertex $v_0$ where $c_G$ has a local extremum is of degree $1$.
\item There exists a piecewise smooth embedding $\tilde{c_G}:G \rightarrow {\mathbb{R}}^2$ satisfying the relation $c_G={\pi}_{2,1} \circ \tilde{c_G}$.
\end{enumerate}
Let $V_G$ be the vertex set of $G$.
Then there exists a Morse refined algebraic domain $({\mathcal{S}}_{\tilde{c_G}},D_{{\mathcal{S}}_{\tilde{c_G}}})$ whose Poincar'e-Reeb V-digraph $(W_{D_{{\mathcal{S}}_{\tilde{c_G}}}},V_{D_{{\mathcal{S}}_{\tilde{c_G}}}})$ for ${\pi}_{2,1}$ is isomorphic to the V-digraph $(G,c_G {\mid}_{V_G})$.
In addition, if the restriction $c_G {\mid}_{V_G}$ is injective, then our Morse refined algebraic domain $({\mathcal{S}}_{\tilde{c_G}},D_{{\mathcal{S}}_{\tilde{c_G}}})$ can be chosen as a pair with ${\mathcal{S}}_{\tilde{c_G}}$ being a one-element set: this has been essentially shown first in \cite{bodinpopescupampusorea} as a main theorem without saying that the curve ${\mathcal{S}}_{\tilde{c_G}}$ is the zero set of some real polynomial and \cite{kitazawa3} (\cite[Remark 4]{kitazawa3}) explains that the non-singular real algebraic curve is also regarded as the zero set of a suitable real polynomial. Note that refined algebraic domains are first defined in \cite{kitazawa6} formally. \cite{bodinpopescupampusorea} is an original motivating study.
We will explain more precise arguments on this in Example \ref{ex:2}, which shows an example for Main Theorem \ref{mthm:2}, presented just after the present example and just before Example \ref{ex:2}. Example \ref{ex:2} and a certain revised case is also important in Main Theorem \ref{mthm:3}, presented later.
\end{Ex}
\begin{MainThm}
\label{mthm:2}
In Main Theorem \ref{mthm:1}, suppose that each connected curve $S_j \in \mathcal{S}$ has no double tangent line tangent at points of $S_j \bigcap \overline{D_{\mathcal{S}}} \bigcap (F_{D_{\mathcal{S}},1,{\rm pole}} \bigcup F_{D_{\mathcal{S}},2,{\rm pole}})$.
Then we can have another Morse refined algebraic domain $({\mathcal{S}}^{\prime},D_{{\mathcal{S}}^{\prime}})$ enjoying the following.
\begin{enumerate}
\item The new refined algebraic domain $({\mathcal{S}}^{\prime},D_{{\mathcal{S}}^{\prime}})$ is NDTL.
\item The relations $\mathcal{S} \subset {\mathcal{S}}^{\prime}$ and $D_{{\mathcal{S}}^{\prime}} \subset D_{\mathcal{S}}$ hold.
\item The V-digraphs $(W_{D_{\mathcal{S}},i},V_{D_{\mathcal{S}},i})$ and $(W_{D_{{\mathcal{S}}^{\prime}},i},V_{D_{{\mathcal{S}}^{prime}},i})$ are homeomorphic for $i=1,2$.
\end{enumerate}

\end{MainThm}
\begin{proof}
Each point $p$ of each curve $S_j$ must not be points of $F_{D_{\mathcal{S}},1,{\rm pole}} \bigcup F_{D_{\mathcal{S}},2,{\rm pole}}$ if it is in some double tangent line of $S_j$.

We choose a point $p$ of a smooth connected curve $S_j \in \mathcal{S}$ in such a way that $p$ is also in the closure $\overline{D_{\mathcal{S}}}$ and a point in a double tangent line of $S_j$.
We choose such a point contained in exactly one curve $S_j \in \mathcal{S}$. 

We apply technique similar to one in the proof of Main Theorem \ref{mthm:1}. 
We can abuse the notation there again. See also FIGUREs \ref{fig:1} and \ref{fig:2} again. For example, the green colored curve connecting the two green colored points $p_{\rm l} \neq p$ and $p_{\rm r} \neq p$ in the second figure from FIGURE \ref{fig:1} is a curve in the closed curve $\partial D_p$.

Note again that the boundary of an ellipsoid does not contain points contained in any double tangent line of the ellipsoid. Note also that our approximation does not change this property.

This completes (exposition on important ingredient of) our proof.
\end{proof}
\begin{Ex}
\label{ex:2}
A case of Example \ref{ex:1} can be also obtained as an explicit case for Main Theorem \ref{mthm:2}.
This is our main result of \cite{kitazawa6}, which is also presented in \cite{kitazawa8} and in Japanese. We explain this.
We change several arguments slightly for our present paper.

We first change the planar graph $\tilde{c_G}(G)$ slightly.
For each vertex $v=(v_1,v_2)$ of the planar graph $\tilde{c_G}(G)$ where $c_G$ does not have a local extremum, we choose a small open neighborhood $U_v \subset {\mathbb{R}}^2$ and another point $v^{\prime}:=(v_1,v_3)$ sufficiently close to the original $v$. We put a straight segment connecting $v$ and $v^{\prime}$. We change the edges departing from $v$ to ones departing from $v^{\prime}$ in $U_v$: remember that $\tilde{c_G}(G)$ is regarded as the V-digraph by the function ${\pi}_{2,1} \circ \tilde{c_G}$ and its restriction to its vertex set. We do this for each of such vertices $v$ of the planar graph $\tilde{c_G}(G)$.  

We can have a small regular neighborhood of the new planar graph $G^{\prime} \subset {\mathbb{R}}^2$ in the smooth category. By considering approximation respecting the class of $C^k$ ($k \geq 2$) the boundary $C_G$ of the regular neighborhood $D_G$ is regarded as the zero set of a real polynomial, non-singular and connected. We can do as follows. This is based on elementary arguments on approximation mainly. See also the original exposition of \cite{kitazawa6, kitazawa8}.
\begin{enumerate}
\item The restriction ${\pi}_{2,1} {\mid}_{C_G}$ is a Morse function. The restriction ${\pi}_{2,1,2} {\mid}_{C_G}$ is a Morse function.
\item Each singular point of the function ${\pi}_{2,1} {\mid}_{C_G}$ corresponds to either the following and is located as presented there.
\begin{enumerate}
\item A vertex $v_0$ of degree $1$. The corresponding singular point of (the function) is located sufficiently close to $c_G(v_0)$.
\item A connected component of $U_v-(U_v \bigcap G^{\prime})$. The corresponding singular point of (the function) is located sufficiently close to $v_0$.
\end{enumerate}
In addition, this gives a one-to-one correspondence.
\item We can add circles and the disks bounded by the circles as follows and remove the intersections of each disk and the interior ${D_G}^{\circ}$ of $D_G$ and we have a desired case. Note that in the original arguments we do not use circles and disks bounded by them and that we use ellipsoids there.
\begin{enumerate}
\item For each singular point $p_v=(p_{v,1},p_{v,2})$ of the previous Morse function corresponding to a connected component of $U_v-(U_v \bigcap G^{\prime})$, we add a small circle bounding the disk whose center is the singular point $p_v$ of the function and whose boundary contains the point of the form $(v_1,p_{v,2})$.  
\item For each singular point $p_{v,0}=(p_{v,0,1},p_{v,0,2})$ of the previous Morse function corresponding to a vertex $v_0$ of degree $1$, we add a small circle bounding the disk which contains $p_{v_0}$ in its interior. We can also add this in such a way that for the intersection of $D_G$ and the circle of the boundary, the restrictions of ${\pi}_{2,1}$ and ${\pi}_{2,1,2}$ there are injective. 
\item We can also choose the disks so that they are mutually disjoint. 
\end{enumerate}
\end{enumerate}

\end{Ex}
In addition, we consider {\it curvatures} for smooth curves in ${\mathbb{R}}^2$.
We can define the {\it curvature} ({\it function}) for a smooth oriented curve (immersion) uniquely as a smooth function.

The curvature function $c:S \rightarrow {\mathbb{R}}^2$ of a smooth oriented curve (immersion) is well-defined as follows. First we consider a local coordinate around $p \in S$ and regard the curve as a smooth map $c_{I_p}=(c_{I_p,1},c_{I_p,2}):I_p \ni p \rightarrow {\mathbb{R}}^2$ on an open interval $I_p$ (which is compatible with the orientation). We consider the $2 \times 2$ matrix whose $(i,j)$-component is the value of the $i$-th derivative of ${c_{I_p,j}}^{(i)}$ at $p$ and let $\det({c_{I_p}}^{(1)},{c_{I_p}}^{(2)})$ denote the determinant of this matrix. The value of our function is $\frac{\det({c_{I_p}}^{(1)},{c_{I_p}}^{(2)})}
{{({({c_{I_p,1}}^{(1)})}^2+{({c_{I_p,1}}^{(1)})}^2)}^{\frac{3}{2}}}$ at $p$.

Thanks to classical theory on differential geometry of smooth curves, the curvature of a smooth curve is constant if and only if it is a subset of a straight line or a subset of a circle: in the former case the values of the function are always $0$ and in the latter case the values of the function are always a fixed value which is not $0$.
A point $p$ of a smooth curve $S$ is a {\it curvature vertex} or a {\it CV of $S$} if it is a singular point of the curvature function for $S$ and an isolated point on a small open neighborhood $U_p$ of $p$. Our arguments here do not depend on the orientation of the curve $S$.
\begin{Def}
If for a refined algebraic domain $(\mathcal{S},D_{\mathcal{S}})$, no CV of a curve $S_j \in \mathcal{S}$ being also a point of $\overline{D_{\mathcal{S}}}-D_{\mathcal{S}}$ exists for any $S_j \in \mathcal{S}$, then the refined algebraic domain is {\it with no CV} or {\it NCV}.
\end{Def}

\begin{MainThm}
\label{mthm:3}
In Main Theorem \ref{mthm:1}, suppose that each connected curve $S_j \in \mathcal{S}$ has no CV which is also a point of $S_j \bigcap \overline{D_{\mathcal{S}}} \bigcap (F_{D_{\mathcal{S}},1,{\rm pole}} \bigcup F_{D_{\mathcal{S}},2,{\rm pole}})$.
Then we can have another Morse refined algebraic domain $({\mathcal{S}}^{\prime},D_{{\mathcal{S}}^{\prime}})$ enjoying the following.
\begin{enumerate}
\item The new refined algebraic domain $({\mathcal{S}}^{\prime},D_{{\mathcal{S}}^{\prime}})$ is NCV.
\item The relations $\mathcal{S} \subset {\mathcal{S}}^{\prime}$ and $D_{{\mathcal{S}}^{\prime}} \subset D_{\mathcal{S}}$ hold.
\item The V-digraphs $(W_{D_{\mathcal{S}},i},V_{D_{\mathcal{S}},i})$ and $(W_{D_{{\mathcal{S}}^{\prime}},i},V_{D_{{\mathcal{S}}^{prime}},i})$ are homeomorphic for $i=1,2$.
\end{enumerate}
\end{MainThm}
\begin{proof}
We can show similalry. We can abuse the notation same as that of the proof of Main Theorems \ref{mthm:1} and \ref{mthm:2}. Most important is that the curve in the given region $\overline{D_{\mathcal{S}}}$ connecting the points $p_{\rm l}, p_{\rm r} \in \overline{D_{\mathcal{S}}}-{D_{\mathcal{S}}}$ can be and must be chosen as a curve with no CV instead of (a curve in the boundary of) an ellipsoid (and a disk $D_p$ obtained by the approximation).

This completes the proof.
\end{proof}
We can add small circles bounded by disks for each point of $F_{D_{\mathcal{S}},2,{\rm pole}}$ similarly, to the obtained example of Example \ref{ex:2}. This shows an example for Main Theorem \ref{mthm:3}.

Last we present related problems abstractly.
\begin{Prob}
We consider several explicit cases of characteristic finite sets of refined algebraic domains respecting differential geometry in our paper. We have also shown related fundamental theorems: Main Theorems \ref{mthm:1}--\ref{mthm:3}. Can we generalize them to a certain generalized theorem?
\end{Prob}
\begin{Prob}
What are meanings of our characteristic finite sets of refined algebraic domains respecting differential geometry in our present study? What are their meanings in real algebraic geometry and singularity theory, for example?
\end{Prob}
\section{Conflict of interest and data availability.}
\noindent {\bf Conflict of interest.}
The author has worked at Institute of Mathematics for Industry (https://www.jgmi.kyushu-u.ac.jp/en/about/young-mentors/). Our present work is closely related to our study. We thank them for their supports and encouragement. The author is also a researcher at Osaka Central
Advanced Mathematical Institute (OCAMI researcher), supported by MEXT Promotion of Distinctive Joint Research Center Program JPMXP0723833165, where he is not employed there. We also thank them for the hospitality. 
The author would also like to thank "Singularity theory of differentiable maps and its applications" (https://www.fit.ac.jp/$\sim$fukunaga/conf/sing202412.html) for letting the author to present our previous result \cite{kitazawa1, kitazawa2, kitazawa3}. This conference is supported by the Research Institute for Mathematical Sciences, an International Joint Usage/Research Center located in Kyoto University and related non-refereed article \cite{kitazawa8} will be published. \\
\ \\
{\bf Data availability.} \\
We do not have data related to the present study other than the present file.


\begin{thebibliography}{25}
%	\bibitem{buchstaberpanov} V. M. Buchstaber and T. E. Panov, \textsl{Toric topology}, Mathematical Surveys and Monographs, Vol. 204, American Mathematical Society, Providence, RI, 2015.
%	\bibitem{burletderham} O. Burlet and G. de Rham, \textsl{Sur certaines applications g\'en\'eriques d'une vari\'et\'e close a $3$ dimensions dans le plan}, Enseign. Math. 20 (1974). 275--292.
	%	\bibitem{calabi} E. Calabi, Quasi-surjective mappings and a generation of Morse theory, Proc. U.S.-Japan Seminar in Differential Geometry, Kyoto, 1965, pp. 13--16.
	%
	%		\bibitem{cavicchioli} A. Cavicchioli, \textsl{Covering numbers of manifolds and critical points of a Morse function}, Israel. J. Math. 70 (1990), 279--304.
	% \bibitem{cerf} J. Cerf, \textsl{La stratification naturelle des espaces de fonctions deff\'erentiables r'eelles et le th'eor`eme de la pseudo-isotopie}, Inst. Hautes Etudes Sci. Publ. Math. 39 (1970), 5--173.
	%		\bibitem{choimasudasuh} S. Choi, M. Masuda and D. Y. Suh, \textsl{Topological classification of generalized Bott towers}, Trans. Amer. Math. Soc. 362 (2010), 1097--1112.
	%		\bibitem{cornealuptonopreatanre} O. Cornea, G. Lupton, J. Oprea and D. Tanr\'e, \textsl{Lusternik-Schnirelmann category}, Mathematical Surveys and Monographs, 103, Amer. Math. Soc., Providence, RI, 2003.
	%\bibitem{crowleyescher} D. Crowley and C. Escher, \textsl{A classification of $S^3$-bundles over $S^4$}, Differential. Geom. Appl. 18 (2003), 363--380, arXiv:0004147.
	%\bibitem{crowleynordstrom} D. Crowley and J. Nordstr\"{o}m, \textsl{The classification of $2$-connected $7$-manifolds}, Proc. London. Math. Soc. 119 (2019), 1--54, arXiv:1406.2226.

	\bibitem{bochnakcosteroy} J. Bochnak, M. Coste and M.-F. Roy, \textsl{Real algebraic geometry}, Ergebnisse der Mathematik und ihrer Grenzgebiete (3) [Results in Mathematics and Related Areas (3)], vol. 36, Springer-Verlag, Berlin, 1998. Translated from the 1987 French original; Revised by the authors.
%		\bibitem{bochnakkucharz} J. Bochnak and W. Kucharz, \textsl{Algebraic approximation of mappings into spheres}, Michigan Mathematical Journal, vol. 34, no. 1, 1987.
	\bibitem{bodinpopescupampusorea} A. Bodin, P. Popescu-Pampu and M. S. Sorea, \textsl{Poincar\'e-Reeb graphs of real algebraic domains}, Revista Matem\'atica Complutense, https://link.springer.com/article/10.1007/s13163-023-00469-y, 2023, arXiv:2207.06871v2.
%\bibitem{bott} R. Bott, \textsl{Nondegenerate critical manifolds}, Ann. of Math. 60 (1954), 248--261.
%\bibitem{carmesinschulz} S. Carmesin and A. Schulz, \textsl{Arrangements of orthogonal circles with many intersections}, Graph Drawing and Network Visualization (a conference paper), SPRINGER NATURE Link, 2021, arXiv:2106.03557v2. 
%\bibitem{costantino}  F. Costantino, \textsl{A short introduction to shadows of $4$-manifolds}, Fundamenta Mathematicae 251 no. 2 (2005), 427--442.
%\bibitem{costantinothurston} F. Costantino, D. Thurston, \textsl{$3$-manifolds efficiently bound $4$-manifolds}, J. Topol. 1 (2008),
%703--745.
%	\bibitem{delzant} T. Delzant, \textsl{Hamiltoniens p\'eriodiques et images convexes de l'application moment}, Bull. Soc. Math. France 116 (1988), No. 3, 315--339.
%\bibitem{ehresmann} C. Ehresmann, \textsl{Les connexions infinitesimales dans un espace fibre differentiable}, Colloque de Topologie, Bruxelles (1950), 29--55.
\bibitem{elredge} N. Elredge, \textsl{{\it Answer to} On finding polynomials that approximate a function and its derivative}, StackExchange, question 555712 (2013), https://math.stackexchange.com/questions/555712/on-finding-polynomials-that-approximate-a-function-and-its-derivative-extension.
%\bibitem{fujitakitabeppumitsuishi} H. Fujita, Y Kitabeppu and A. Mitsuishi, \textsl{Distance functions and convex bodies and symplectic toric manifolds}, arXiv:2003.02293.
%\bibitem{gelbukh} I. Gelbukh, \textsl{Loops in Reeb graphs of $n$-manifolds}, diskrete \& Computational Geometry, 59 (4) (2018), 843--863. 
%%\bibitem{gelbukh2} I. Gelbukh, \textsl{Approximation of Metric Spaces by Reeb Graphs: Cycle Rank of a Reeb Graph, the Co-rank of the Fundamental Group, and Large Components of Level Sets on Riemannian Manifolds}, Filomat (in press), arxiv:1903.00777.
%\bibitem{gelbukh1} I. Gelbukh, \textsl{A finite graph is homeomorphic to the Reeb graph of a Morse-Bott function}, Mathematica Slovaca, 71 (3), 757--772, 2021; doi: 10.1515/ms-2021-0018. 
%\bibitem{gelbukh2} I. Gelbukh, \textsl{Morse-Bott functions with two critical values on a surface}, Czechoslovak Mathematical Journal, 71 (3), 865--880, 2021; doi: 10.21136/CMJ.2021.0125-20. 
\bibitem{golubitskyguillemin} M. Golubitsky and V. Guillemin, \textsl{Stable Mappings and Their Singularities}, Graduate Texts in Mathematics (14), Springer-Verlag (1974).
%\bibitem{hempel} J. Hempel, \textsl{3- Manifolds}, AMS Chelsea Publishing, 2004. 
%\bibitem{hiratukasaeki} J. T. Hiratuka and O. Saeki, \textsl{Triangulating Stein factorizations of generic maps and Euler Characteristic formulas}, RIMS Kokyuroku Bessatsu B38 (2013), 61--89. 
%\bibitem{hiratukasaeki2} J. T. Hiratuka and O. Saeki, \textsl{Connected components of regular fibers of differentiable maps}, in "Topics on Real and Complex Singularities", Proceedings of the 4th Japanese-Australian Workshop (JARCS4), Kobe 2011,  World Scientific, 2014, 61--73. 
%\bibitem{hirsch} M. W. Hirsch, \textsl{Smooth regular neighborhoods}, Ann. of Math., 76 (1962), 524--530.
%\bibitem{ishikawakoda} M. Ishikawa and Y. Koda, \textsl{Stable maps and branched shadows of $3$-manifolds}, Mathematische Annalen 367 (2017), no. 3, 1819--1863, arXiv:1403.0596.
%\bibitem{kitazawa1} N. Kitazawa, \textsl{On round fold maps} (in Japanese), RIMS Kokyuroku Bessatsu B38 (2013), 45--59.
%\bibitem{kitazawa2} N. Kitazawa, \textsl{On manifolds admitting fold maps with singular value sets of concentric spheres}, Doctoral Dissertation, Tokyo Institute of Technology (2014).
%\bibitem{kitazawa3} N. Kitazawa, \textsl{Fold maps with singular value sets of concentric spheres}, Hokkaido Mathematical Journal Vol.43, No.3 (2014), 327--359.
%\bibitem{kitazawa1} N. Kitazawa, \textsl{On Reeb graphs induced from smooth functions on $3$-dimensional closed orientable manifolds with finitely many singular values}, Topol. Methods in Nonlinear Anal. Vol. 59 No. 2B, 897--912, arXiv:1902.08841.
%\bibitem{kitazawa1} N. Kitazawa, \textsl{On Reeb graphs induced from smooth functions on closed or open surfaces}, Methods of Functional Analysis and Topology Vol. 28 No. 2 (2022), 127--143, arXiv:1908.04340.
\bibitem{kitazawa1} N. Kitazawa, \textsl{Real algebraic functions on closed manifolds whose Reeb graphs are given graphs}, Methods of Functional Analysis and Topology Vol. 28 No. 4 (2022), 302--308, arXiv:2302.02339, 2023.
%\bibitem{kitazawa6} N. Kitazawa, \textsl{Explicit construction of explicit real algebraic functions and real algebraic manifolds via Reeb graphs}, Algebraic and geometric methods of analysis 2023 “The book of abstracts”, 49—51, this is the abstract book of the conference "Algebraic and geometric methods of analysis 2023" and published after a short review (https://www.imath.kiev.ua/$\sim$topology/conf/agma2023/), https://imath.kiev.ua/$\sim$topology/conf/agma2023/contents/abstracts/texts/kitazawa/kitazawa.pdf, 2023.
%\bibitem{kitazawa5} N. Kitazawa, \textsl{Notes on explicit special generic maps into Euclidean spaces whose dimensions are greater than $4$}, a revised version is submitted based on positive comments (major revision) by referees and editors after the first submission to a refereed journal, arXiv:2010.10078.

%\bibitem{kitazawa6} N. Kitazawa, \textsl{Round fold maps and the topologies and the differentiable structures of manifolds admitting explicit ones}, submitted to a refereed journal, arXiv:1304.0618.
%\bibitem{kitazawa0.5} N. Kitazawa, \textsl{Constructing fold maps by surgery operations and homological information of their Reeb spaces}, submitted to a refereed journal, arxiv:1508.05630.
%\bibitem{kitazawa0.6} N. Kitazawa, \textsl{Notes on fold maps obtained by surgery operations and algebraic information of their Reeb spaces}, arxiv:1811.04080.


%\bibitem{kitazawa6} N. Kitazawa, \textsl{On Reeb graphs induced from smooth functions on $3$-dimensional closed manifolds which may not be orientable}, a revised version is submitted to a refereed journal after based on positive comments by editors and referees after the second submission to a refreed journal, arXiv:2108.01300.
%\bibitem{kitazawa7} N. Kitazawa, \textsl{Realization problems of graphs as Reeb graphs of Morse functions with prescribed preimages}, submitted to a refereed journal, arXiv:2108.06913.
%\bibitem{kitazawa10} N. Kitazawa,\textsl{Round fold maps on $3$-dimensional manifolds and their integral and rational cohomology rings}, arXiv:2301.07008.
%\bibitem{kitazawa6} N. Kitazawa, \textsl{A class of naturally generalized special generic maps}, arXiv:2212.03174.
%\bibitem{kitazawa7} N. Kitazawa, \textsl{Construction of real algebraic functions with prescribed preimages}, submitted to a refereed journal as the second version based on positive comments by referees and editors, arXiv:2303.00953v3.
\bibitem{kitazawa2} N. Kitazawa, \textsl{Reconstructing real algebraic maps locally like moment-maps with prescribed images and compositions with the canonical projections to the $1$-dimensional real affine space}, the title has changed from previous versions, arXiv:2303.10723, 2024.
%\bibitem{kitazawa3} N. Kitazawa, \textsl{Explicit smooth real algebraic functions which may have both compact and non-compact preimages on smooth real algebraic manifolds}, arXiv:2304.07450.

%\bibitem{kitazawa6} N. Kitazawa, \textsl{A note on real algebraic maps which are topologically special generic maps}, previous version(s) of the present article and the version arXiv:2303.00953v2 is submitted to a refereed journal, arXiv:2312.10646. 
\bibitem{kitazawa3} N. Kitazawa, \textsl{Some remarks on real algebraic maps which are topologically special generic maps}, arXiv:2312.10646. 
\bibitem{kitazawa4} N. Kitazawa, \textsl{Arrangements of small circles for Morse-Bott functions and regions surrounded by them}, arXiv:2412.20626v3, 2025.
\bibitem{kitazawa5} N. Kitazawa, \textsl{Arrangements of circles supported by small chords and compatible with natural real algebraic functions}, arXiv:2501.11819.
\bibitem{kitazawa6} N. Kitazawa, \textsl{Realizations of planar graphs of refined algebraic domains}, arXiv:2501,17425, 2025.
\bibitem{kitazawa7} N. Kitazawa, \textsl{Refined algebraic domains with finite sets in the boundaries}, arXiv:2503.09195 ,2025. 
\bibitem{kitazawa8} N. Kitazawa, \textsl{Constructing real algebraic functions explicitly and their Reeb graphs}, an non-refereed article on the talk of the same title in the conference "Singularity theory of differentiable maps and its applications" (https://www.fit.ac.jp/$\sim$fukunaga/conf/sing202412.html) and submitted to RIMS K\^oky\^uroku.
%\bibitem{kitazawa10} N. Kitazawa, \textsl{A note on cohomological structures of special generic maps}, a revised version is submitted based on positive comments by referees and editors after the third submission to a refereed journal.
%\bibitem{kitazawasaeki1} N. Kitazawa and O. Saeki, \textsl{Round fold maps on $3$-manifolds}, accepted for publication after a refereeing process and to appear in Algebraic \& Geometric Topology, arXiv:2105.00974.
		%	\bibitem{kitazawasaeki2} N. Kitazawa and O. Saeki, \textsl{Round fold maps of $n$-dimensional manifolds into ${\mathbb{R}}^{n-1}$}, submitted to a refereed journal, arXiv:2111.13510.
%\bibitem{ishikawakoda} M. Ishikawa, Y. Koda, \textsl{Stable maps and branched shadows of $3$-manifolds}, arXiv:1403.0596.
%\bibitem{kobayashisaeki} M. Kobayashi and O. Saeki, \textsl{Simplifying stable mappings into the plane from a global viewpoint}, Trans. Amer. Math. Soc. 348 (1996), 2607--2636.
%\bibitem{kohnpieneranestadrydellshapirosinnsoreatelen} K. Kohn, R. Piene, K. Ranestad, F. Rydell, B. Shapiro, R. Sinn, M-S. Sorea and S. Telen, \textsl{Adjoints and Canonical Forms of Polypols}, to appear in Documenta Mathematica, arXiv:2108.11747.
\bibitem{kollar} J. Koll\'ar, \textsl{Nash's work in algebraic geometry}, Bulletin (New Series) of the American Mathematical Society (2) 54, 2017, 307--324.
%\bibitem{kucharz} W. Kucharz, \textsl{Some open questions in real algebraic geometry}, Proyecciones Journal of Mathematics, Vol. 41 No. 2 (2022), Universidad Cat\'olica del Norte Antofagasta, Chile, 437--448.
\bibitem{lellis} Camillo De Lellis, \textsl{The Masterpieces of John Forbes Nash Jr.}, The Abel Prize 2013--2017 (Helge Holden and Ragni Piene, eds.), Springer International Publishing, Cham, 2019, 391--499, https://www.math.ias.edu/delellis/sites/math.ias.edu.delellis/files/Nash\_Abel\_75.pdf, arXiv:1606.02551.
%\bibitem{martinezalfaromezasarmientooliveira} J. Martinez-Alfaro, I. S. Meza-Sarmiento and R. Oliveira, \textsl{Topological classification of simple Morse Bott functions on surfaces}, Contemp. Math. 675 (2016), 165--179.%
%\bibitem{marzantowiczmichalak} W. Marzantowicz and L. P. Michalak, \textsl{Relations between Reeb graphs, systems of hypersurfaces and epimorphisms onto free groups}, Fund. Math., 265 (2), 97--140, 2024.
%\bibitem{masumotosaeki} Y. Masumoto and O. Saeki, \textsl{A smooth function on a manifold with given Reeb graph}, Kyushu J. Math. 65 (2011), 75--84.
%\bibitem{maciasvirgospereirasaez} E. Mac\'ias-Virg\'os and M. J. Pereira-S\'aez, Height functions on compact symmetric spaces, Monatshefte f\"ur Mathematik 177 (2015), 119--140. 
%\bibitem{michalak1} L. P. Michalak, \textsl{Realization of a graph as the Reeb graph of a Morse function on a manifold}. Topol. Methods in Nonlinear Anal. 52 (2) (2018), 749--762, arXiv:1805.06727.
%\bibitem{michalak2} L. P. Michalak, \textsl{Combinatorial modifications of Reeb graphs and the realization problem}, Discrete Comput. Geom. 65 (2021), 1038--1060, arXiv:1811.08031.
%\bibitem{milnor} J. Milnor, \textsl{Singular points of complex hypersurfacs}, Annals of Mathematics Studies, No. 61, Princeton University Press, Princeton, N. J.; University of Tokyo Press, Tokyo, 1968.
%\bibitem{milnor} J. Milnor, \textsl{Lectures on the h-cobordism theorem}, Math. Notes, Princeton Univ. Press, Princeton, N.J. 1965.
%\bibitem{moise} E. E. Moise, \textsl{Affine Structures in $3$-Manifold{\rm :} V. The Triangulation Theorem and Hauptvermutung}, Ann. of Math., Second Series, Vol. 56, No. 1 (1952), 96--114.
%\bibitem{morin} B. Morin, \textsl{Formes canoniques des singulariti\'{e}s d\'{}une application diff\'{e}rentiable}, C. E. Acad. Sci. Paris 260 (1965), 5662--5665, 6503--6506.
%\bibitem{nash} J. Nash, \textsl{Real algbraic manifolds}, Ann. of Math. (2) 56 (1952), 405--421.
%\bibitem{ranicki} A. Ranicki, \textsl{Algebraic and geometric surgery}, https://www.maths.ed.ac.uk/~v1ranick/books/surgery.pdf, 2002.
%\bibitem{ramanujam} S. Ramanujam, \textsl{Morse theory of certain symmetric spaces}, J. Diff. Geom. 3 (1969), 213--229.
%\bibitem{reeb} G. Reeb, \textsl{Sur les points singuliers d\'{}une forme de Pfaff compl\'{e}tement int\`{e}grable ou d\'{}une fonction num\'{e}rique}, Comptes Rendus
% Hebdomadaires des S\'{e}ances de I\'{}Acad\'{e}mie des Sciences 222 (1946), 847--849.
%\bibitem{saeki1} O. Saeki, \textsl{Notes on the topology of folds}, J. Math. Soc. Japan Volume 44, Number 3 (1992), 551--566.
%\bibitem{saeki1} O. Saeki, \textsl{Topology of special generic maps of manifolds into Euclidean spaces}, Topology Appl. 49 (1993), 265--293.
%\bibitem{saeki0.2} O. Saeki, \textsl{Topology of singular fibers of differentiable maps}, Lecture Notes in Math., Vol. 1854, Springer-Verlag, 2004. 
%\bibitem{saeki4} O. Saeki, \textsl{Morse functions with sphere fibers}, Hiroshima Math. J. Volume 36, Number 1 (2006),  141--170.
\bibitem{saeki1} O. Saeki, \textsl{Reeb spaces of smooth functions on manifolds}, International Mathematics Research Notices, maa301, Volume 2022, Issue 11, June 2022, 3740--3768, https://doi.org/10.1093/imrn/maa301.

\bibitem{saeki2} O. Saeki, \textsl{Reeb spaces of smooth functions on manifolds II}, Res. Math. Sci. 11, article number 24 (2024), https://link.springer.com/article/10.1007/s40687-024-00436-z.
%\bibitem{saekitakase} O. Saeki and M. Takase, \textsl{Desingularizing special generic maps}, Journal of G\"okova Geometry Topology (2013), 1--24.
%\bibitem{sakurai} S. Sakurai, Master Thesis, Kyushu. Univ..
% \bibitem{saekitakase} O. Saeki and M. Takase, \textsl{Desingularizing special generic maps}, Journal of Gokova Geometry Topology 7 (2013), 1--24.
%\bibitem{saeki2} O. Saeki, \textsl{Topology of special generic maps of manifolds into Euclidean spaces}, Topology Appl. 49 (1993), 265--293.
%\bibitem{saeki4} O. Saeki, \textsl{Singular fibers and $4$-dimensional cobordism group}, Pacific J. Math. 248 (2010), 233--256.
%\bibitem{saekisakuma} O. Saeki and K. Sakuma, \textsl{On special generic maps into ${\mathbb{R}}^3$}, Pacific J. Math. 184 (1998), 175--193.
%\bibitem{saekisuzuoka} O. Saeki and K. Suzuoka, \textsl{Generic smooth maps with sphere fibers} J. Math. Soc. Japan Volume 57, Number 3 (2005), 881--902.
%\bibitem{sharko} V. Sharko, \textsl{About Kronrod-Reeb graph of a function on a manifold}, Methods of Functional Analysis and
 Topology 12 (2006), 389--396.
%\bibitem{shiota} M. Shiota, \textsl{Thom's conjecture on triangulations of maps}, Topology 39 (2000), 383--399.
\bibitem{sorea1} M. S. Sorea, \textsl{The shapes of level curves of real polynomials near strict local maxima},  Ph. D. Thesis, Universit\'e de Lille, Laboratoire Paul Painlev\'e, 2018.
\bibitem{sorea2} M. S. Sorea, \textsl{Measuring the local non-convexity of real algebraic curves}, Journal of Symbolic Computation 109 (2022), 482--509.
%\bibitem{sorea1} M. S. Sorea, \textsl{The shapes of level curves of real polynomials near strict local maxima},  Ph. D. Thesis, Universit\'e de %Lille, Laboratoire Paul Painlev\'e, 2018.
%\bibitem{sorea2} M. S. Sorea, \textsl{Measuring the local non-convexity of real algebraic curves}, J. Symbolic Compute. 109 (2022), 482--509.
%\bibitem{stong} R. E. Stong, \textsl{Notes on cobordsm theory}, Princeton Universty Press, 1968.
%\bibitem{takeuchi} M. Takeuchi, \textsl{Nice functions on symmetric spaces}, Osaka. J. Mat. (2) Vol. 6 (1969), 283--289.
%\bibitem{tamaki1} D. Tamaki, Algebraic Topology A Guide to literature,  http://pantodon.jp/index.rb?body=about, 2023. 
%\bibitem{tamaki2} D. Tamaki, Algebraic Topology A Guide to literature (Submanifold arrangement), http://pantodon.jp/index.rb?body=submanifold\_arrangement, 2023.
%\bibitem{tamaki3} D. Tamaki, Algebraic Topology A Guide to literature (Arrangement variations), http://pantodon.jp/index.rb?body=arrangement\_variations, 2023.

%\bibitem{thom} R. Thom, \textsl{Les singularites des applications differentiables}, Ann. Inst. Fourier (Grenoble) 6 (1955-56), 43--87.
%\bibitem{tognoli} A. Tognoli, \textsl{Su una congettura di Nash}, Ann. Scuola Norm. Sup. Pisa (3) 27 (1973), 167--185.
%\bibitem{turaev} Vladimir G. Turaev, \textsl{Topology of shadows}, Preprint, 1991.
%\bibitem{wall} C. T. C Wall, \textsl{Classification problems in differential topology -- {\rm I:} Classificationon handlebodies}, Topology 2 (1963), 253--261.
%\bibitem{wall2} C. T. C. Wall \textsl{Classification problems in differential topology -- {\rm II:} Diffeomorphismsof handlebodies}, Topology 2 (1963), 263--272.
%\bibitem{wall3} C. T. C. Wall, \textsl{Classification problems in differential topology -- {\rm Q:} Quadratic forms on finite groups and related topics}, Topology 2 (1963), 281--298.
%\bibitem{wall4} C. T. C. Wall, \textsl{Classification problems in differential topology -- {\rm III:} Applications to special cases}, Topology 3 (1965), 291--304.
%%\bibitem{wall5} C. T. C. Wall, \textsl{Classification problems in differential topology -- {\rm IV:} Thickenings}, Topology 5 (1966), 73--94.
%\bibitem{wall6} C. T. C. Wall, \textsl{Classification problems in differential topology -- {\rm VI:} Classification of |{\rm (}$s-1${\rm )}-connected {\rm (}$2s+1${\rm )}-manifolds}, Topology 6 (3) (1967), 273--296.
%\bibitem{whitney} H.  Whitney,  \textsl{On singularities of mappings of Euclidean spaces: I,  mappings of the plane into the plane},  Ann.  of Math.  62 (1955),  374--410. 

	%\bibitem{zhubr1} A. V. Zhubr, Closed simply-connected six-dimensional manifolds: proofs of classification theorems, Algebra i Analiz 12 (2000), no. 4, 126--230.
%\bibitem{zhubr2} A. V. Zhubr (responsible for the page), http://www.map.mpim-bonn.mpg.de/6-manifolds:\_1-connected.
\end{thebibliography}
\end{document}